\documentclass[12pt]{amsart}

 \usepackage{amsmath}
\usepackage{amssymb,latexsym,cite}
\usepackage{amscd}
\usepackage{comment}
 \usepackage{xspace}
 \usepackage{verbatim}
 \usepackage{todonotes}

\newtheorem{theorem}{Theorem}[section]
\newtheorem{proposition}[theorem]{Proposition}

\newtheorem{lemma}[theorem]{Lemma}
\theoremstyle{definition}

\numberwithin{equation}{section}

\newcommand{\R}{{\mathbb R}}
\newcommand{\E}{{\mathcal E}}

\newcommand{\G}{\mathbb{G}}

\newcommand{\tr}{\mathrm{tr}^*}

\newcommand{\chr}[1]{\mathbf{1}\ind{#1}}

\newcommand{\BSA}{\begin{subarray}}
\newcommand{\ESA}{\end{subarray}}
\newcommand{\BAL}{\begin{aligned}}
\newcommand{\EAL}{\end{aligned}}

\newcommand{\forevery}{\quad \forall}

\newcommand{\norm}[1]{\left \|#1\right \|}

\newcommand{\rec}[1]{\frac{1}{#1}}

\newcommand{\supp}{\mathrm{supp}\,}
\newcommand{\dist}{\mathrm{dist}\,}

\newcommand{\diam}{\mathrm{diam}\,}

\newcommand{\sms}{\setminus}

\newcommand{\ti}{\times}

\newcommand{\sbs}{\subset}

\newcommand{\ind}[1]{_{_{#1}}}

\newcommand{\sth}{such that\xspace}

\newcommand{\bdw}{\partial\Gw}

\newcommand{\qtxt}[1]{\quad\textrm{#1}}

\def\ga{\alpha}     \def\gb{\beta}       \def\gg{\gamma}
\def\gc{\chi}       \def\gd{\delta}

      \def\gk{\kappa}      \def\gl{\lambda}
                 
    \def\gr{\rho}        
\def\gs{\sigma}       
      
\def\gx{\xi}                
     \def\Gd{\Delta}      

    \def\Gs{\Sigma}      
\def\Gw{\Omega}              

      \def\CC{{\mathcal C}}

\def\BBG {\mathbb G}

   \def\BBR {\mathbb R}

\def\GTM {\mathfrak M}

\def\tr{\mathrm{tr}}

\def\hr{harmonic\xspace}
\def\suphr{superharmonic\xspace}

\def\LVsup{$L_V$ superharmonic\xspace}
\def\LVsub{$L_V$ subharmonic\xspace}
\def\LVhr{$L_V$ harmonic\xspace}

\def\q{\quad}
\def\Lip{Lipschitz\xspace}
\def\mbfn{\mathbf{n}}

\begin{document}

\title[Schr\"odinger equations]{Schr\"odinger equations with very singular potentials in Lipschitz domains.}

\author{Moshe Marcus }
\address{Department of Mathematics, Technion\\
	Haifa 32000, ISRAEL}
\email{marcusm@math.technion.ac.il}

\date{}

\begin{abstract}
	Consider operators  $L^{V}:=\Gd + V$ in a bounded Lipschitz domain $\Gw\sbs \mathbb{R}^N$. Assume that $V\in C^{1,1}(\Gw)$ and $V$ satisfies $V(x) \leq \bar a\,\dist(x,\bdw)^{-2}$ in $\Gw$ and also a second condition that guarantees the existence of a ground state $\Phi_V$. If, for example, $V>0$ this condition reads $1<c_H(V)$ (= the Hardy constant relative to $V$). We derive estimates of positive $L_V$ harmonic functions and of positive Green potentials of measures $\tau\in \GTM_+(\Gw;\Phi_V)$. These imply estimates of positive \LVsup functions and of \LVsub functions that are dominated by an \LVsup. Similar results have been obtained in \cite{MM2019N} in the case of smooth domains.   \\ [2mm]
	
	MSC:35J60; 35J75; 35J10\\ [2mm]
	
	Keywords: $L_V$ potentials, Green and Martin kernels, Boundary Harnack princiole.

\end{abstract}

\dedicatory{This paper is dedicated to Profesor Shmuel Agmon with gratitude and admiration.}

\maketitle

\section{Introduction}
Let $\Gw$ be a bounded, Lipschitz domain in $\R^N$, $N\geq3$. Let
$$L_V:=\Gd +V$$
where $V \in C^{1,1}(\Gw)$.  
We assume that the
potential $V$ satisfies the conditions:
\begin{equation}\label{Vcon1}
\exists \bar a>0\, :\quad |V(x)| \leq \bar a \gd(x)^{-2} \forevery x\in \Gw \tag{A1}
\end{equation}
$$ \gd(x)=\gd_\Gw(x):=\dist(x,\bdw).$$
and,
\begin{equation}\label{Vcon2}
\gg_-<1<\gg_+. \tag{A2}
\end{equation}
where,
\begin{equation}\label{gg+-}\BAL
\gg_+&:=\sup  \{\gg: \exists u_\gg>0\;\text{\sth}\; L^{\gg V}u_\gg=0\},\\
\gg_-&:=\inf\,  \{\gg: \exists u_\gg>0\;\text{\sth}\; L^{\gg V}u_\gg=0\}.
\EAL\end{equation}

By a theorem of Allegretto and Piepenbrink (see \cite{Simon}), \eqref{gg+-} is equivalent to,

\begin{equation}\label{ggpm}
\begin{aligned}
\gg_+&=\sup\{\gg: \int_\Gw |\nabla\phi|^2\,dx\geq \gg\int_\Gw \phi^2 V\, dx  \forevery \phi\in H^1_0(\Gw)\},\\
\gg_-&=\inf\,\{\gg: \int_\Gw |\nabla\phi|^2\,dx\geq \gg\int_\Gw \phi^2 V\, dx  \forevery \phi\in H^1_0(\Gw)\}.
\end{aligned}
\end{equation}

Therefore, condition (A1) and Hardy's inequality imply that $\gg_+>0$ and $\gg_-<0$.
If $V$ is positive then $\gg_-=-\infty$ and $\gg_+$ is \emph{the Hardy constant relative to $V$} in $\Gw$, denoted by $c_H(V)$.
If $V$ is negative then $\gg_+=\infty$.

If $\gg\in (\gg_-,\gg_+)$, the operator    $L_{\gg V}$ has a Green function in $\Gw$. In particular condition (A2) implies that $L_V$ has a Green function $G_V$ in $\Gw$.
Moreover, $L_V$ has  a positive eigenfunction $\Phi_V$ corresponding to the smallest eigenvalue $\gl_V$ of $L_V$. (In the present case $\gl_V>0$.) The eigenfunction is  normalized by setting $\Phi_V(x_0)=1$ where $x_0$ is a fixed reference point in $\Gw$.

The family of potentials satisfying conditions (A1) and (A2) contains in particular the following:
\begin{equation}\label{e:V_E}
V= \gg V_E, \q V_E= \rec{\gd_E^2},\q \gd(E)(x)=\dist (x,\E)
\end{equation}
where $E$ is a closed subset of $\bdw$ and $\gg<c_H(V_E)$ (= Hardy constant, relative to $V_E$, in $\Gw$). 

We list here a few terms that are used below.

	(i) A function $u>0$, \LVhr in a one-sided neighborhood of $\bdw$, is called a \emph{local \LVhr}. A \emph{local \LVsup} is similarly defined.
	
	(ii) A positive local \LVsup function $u$, is a \emph{ground state} if it is of minimal growth in the sense of Agmon \cite{Ag}), i.e., for every positive local $L_V$ \suphr function $v$,
	$$\limsup_{\gd(x)\to0}  \frac{u}{v}(x)<\infty.$$
	When (A2) holds $\Phi_V$ is a ground state.	
	
	(iii) A positive \LVsup is called an $L_V$ \emph{potential} if it does not dominate any positive \LVhr function.\\

	It is known \cite{An-SLN} that $u$ is an $L_V$ potential if and only if
	$$u=\BBG_V[\tau]:=\int_\Gw G_V(x,y) d\tau(y)$$
	where $\tau$ is a positive Radon measure  on $\Gw$ \sth $\BBG_V[\tau]<\infty$.

The Green kernel $G_V$ is uniquely determined by the following conditions, \cite{Ag}.

(a) For every $y\in \Gw$, $-L_V G_V(x,y)= \gd_y$ (the Dirac measure at $y$) and

(b) the function $x\mapsto G_V(x,y)$ is of minimal growth in $\Gw\sms\{y\}$.

\medskip

Conditions (A1) and (A2) imply that $L_V$ is weakly coercive in the sense of Ancona \cite{An87}. For a proof, based on \cite{Pincho-comm}, see \cite[Lemma 1.1]{MM-Green}.
(Recall  that, if $V<0$,  $\gg_+=\infty$ while $\gg_-<0$. Thus (A2) always holds in this case.)
Consequently the results of Ancona \cite{An87} apply to the operators $L_V$ under consideration.  In particular: \\ [2mm]
$\bullet$ There exista  a Martin kernel $K_V$ such that, for every $y\in \bdw$, $x\mapsto K_V(x,y)$ is positive $L_V$ harmonic in $\Gw$ and vanishes on $\bdw\sms \{y\}$ and the following holds,\\ [2mm]
\noindent\textbf{Representation Theorem.}\hskip 2mm
\textit{For every positive $L_V$-\hr function $u$ there exists a measure
	$\nu\in \GTM_+(\bdw)$ ($=$ the space of positive, bounded Borel measures) \sth
	\begin{equation}\label{rep}
	u(x)=\int_\Gw K_V(x,y) d\nu(y)=:K_V[\nu] \q x\in \Gw.
	\end{equation}
	Conversely, for every such measure $\nu$, the function $u$ above is $L_V$ harmonic.}\\ [2mm]
$\bullet$ The \textbf{Boundary Harnack Principle} (briefly BHP) holds. (For its statement - in the form used here - see \cite{MM-Green}.)\\ [2mm]
 $\bullet$ There exists a constant $C > 0$ such that for every $y_0\in \Gw$,
\begin{equation}\label{G-Phi}
C^{-1} G_V(x,y_0)\leq \Phi_V(x)\leq C G_V(x,y_0) \qtxt{when }\; \gd(x)<\rec{2}\gd(y_0).
\end{equation}
This is a consequence of Theorem 3.2 and Lemma 3.6 of \cite{Pincho1989}. 

Inequality \eqref{G-Phi} implies that, if $\tau$ is a positive Radon measure in $\Gw$ then
\begin{equation}
G_V[\tau]<\infty \Longleftrightarrow \tau\in \GTM(\Gw;\Phi_V) \qtxt{i.e. }\int_\Gw\Phi_Vd\tau<\infty.
\end{equation}

\vskip 1mm

Using these results,  sharp estimates of Green and Martin kernels of $L_V$ have recently been obtained  in \cite{MM-Green}.

\vskip 3mm

In a recent paper \cite{MM2019N} we derived sharp estimates of $L_1(\Gw;W)$ norms of positive \LVhr functions and $L_V$ potentials in domains $\Gw$ of class $C^2$ and potentials satisfying (A1), and (A2) and some conditions on the behavior of $\Phi_V$ that are satisfied by a large class of potentials. Here $L_1(\Gw;W)$ denotes the weighted $L_1$ space with weight,
$$   W:=\frac{\Phi_V}{\Phi_0}.$$

In the case of Hardy potentials, $V=\mu/\gd^2$, estimates of $L_V$ harmonic functions in smooth domains and related results for semilinear equations  have been obtained in several papers, see e.g.  \cite{MT2}, \cite{MT3}, \cite{Gk-Ve}, \cite{GkNg}, \cite{MM-VM} and the references therein.

In this paper we show that, with some modifications, the estimates  of \cite{MM2019N} are also valid in bounded Lipschitz domains. For Lipschitz domains the estimates are new even in the case of the Hardy potential. As in the previous paper, the proofs are based on potential theoretic results, quoted above, and the estimates of Green and Martin kernels obtained in \cite{MM-Green}. The global estimates are based on local estimates, i.e., for each point $P\in \bdw$ we derive estimates in an one-sided neighborhood of $P$ in terms of a local set of coordinates centered at $P$. Some of the local estimates don't have a natural global extension (except in a more restricted class of domains) but they provide additional information. For instance, the local estimates of a Green potential  $G_V[\tau]$, $\tau\in \GTM(\Gw;\Phi_V)$, show that weighted integrals on  surfaces parallel to the boundary  tend to zero at the boundary at the rate of the decay of $\tau$ at the boundary  in the neighborhood of $P$. Similarly,  the local estimates of an \LVhr function $K_V[\nu]$ show that, near the boundary, the corresponding weighted surface integrals are equivalent to the norm of $\nu$ in the neighborhood of $P$.  

As a consequence of the estimates presented here, we obtain sharp estimates of \LVsup and \LVsub functions. In turn these can be applied in the study of semilinear equations - with absorption or source nonlinear term - involving a strongly singular potential. These applications will be discussed in a subsequent paper.

The main results are stated in section 2. The proofs are provided in sections 3 -- 5.

\section{Preliminaries and main results}\label{S:Lip}
Let $\Gw$ be a bounded \Lip domain. 
In the present case the flow coordinates $(\gd,\gs)$ - that have been extensively used in \cite{MM-Green} - are not available in general. Instead we shall use local coordinates, defined below. Some of our estimates will be expressed in these local coordinates, but the main results are independent of the local coordinates.

Denote,
$$T(r,\gr) : = \{\gx= (\gx_1, \gx')\in \R\ti\R^{N-1}:    |\gx_1|<\gr,\;  |\gx'|<r\}.$$
For a point $P\in \BBR^N$, denote $T^P(r,\gr):= P+ T(r,\gr)$.

Since $\Gw$ is a bounded Lipschitz domain, there exist positive numbers  $r_0$ , $\gk$ \sth, for every $P\in \bdw$, there exist: (i) a  set of Euclidean coordinates $\gx=\gx^P$ centered at $P$ with the positive $\gx_1$ axis pointing into the domain and (ii) a function $F_P$ uniformly Lipschitz in $\R^{N-1}$ with Lipschitz constant $\leq \kappa$ \sth
\begin{equation}\label{ngh_P}\BAL
Q^P(r_0,\gr_0):=\, &\Gw\cap T^P(r_0,\gr_0)\\
=\, &\{ \gx= (\gx_1, \gx'): F_P(\gx')<\gx_1<\gr_0,\; |\gx'|<r_0\},
\EAL\end{equation}
where $\gr_0 = 10\gk r_0$.
Without loss of generality, we assume that $\kappa>1$. The pair $(r_0,\gk)$ is called the \emph{Lipschitz characteristic} of $\Gw$.

The set of coordinates $\gx^P$ is called a standard set of coordinates at $P$ and $T^P(r,\gr)$ with $0<r\leq r_0$ and $\gr=c\gk r$, $2<c\leq 10$ is called a standard cylinder at $P$. For $y\in T^P(r_0,\gr_0)$, the $\gx^P$ coordinates of $y$ are denoted by $\gx^P(y)$. If $y\in 
T^P(r_0,\gr_0)\cap\bdw$ then $\gx^P_1(y)= F_P((\gx^P)'(y))$.

Let $P\in \bdw$, let $\gx$ be a standard set of coordinates at $P$ and let $T_P(r_0,\gr_0)$ be the corresponding  standard cylinder at $P$. Denote
\begin{equation}\label{d:loc-co}\BAL
\Gs^P_\gb(r,\rho)&:= \{\gx\in Q^P(r,\gr): \gx_1=\gb + F_P(\gx')\},\\
Q_\gb^P(r,\rho) &:=  \{\gx\in Q^P(r,\gr): \gx_1<\gb + F_P(\gx')\},\\
\Gw^P_\gb(r,\rho)&:=  \{\xi\in Q^P(r,\gr): \gd(\xi)<\gb\},\q \gd(\xi)=\dist (\xi,\bdw),
\EAL
\end{equation}
with $F_P$ as in \eqref{ngh_P}, $r\in (0,r_0)$ and $\rho\in (2\gk r, 10\gk r)$. 
Finally denote, $$dS^P_\gb:=dH_{N-1} \qtxt{on }\; \Gs_\gb^P.$$

It is easily verified that there exists a constant $\bar c=\bar c(\kappa)$ \sth, for every $P\in \bdw$ and every $\gb\in (0,\gr_0/2)$,
\begin{equation}\label{Gsy}
\rec{\bar c_\gk}\gb\leq \gd(\gx) \leq  \gb \forevery \gx\in \Gs_\gb^P(r_0/2,\gr_0) \q \text{where }\; \bar c_\gk= \sqrt{1+\gk^2}.
\end{equation}

Let
$$\CC= \{\xi=(\xi_1,\xi')\in \BBR\ti\BBR^{N-1}:  |\xi'|<\gk \xi_1, 0<\xi_1<r_0/\gk \}.$$ 
If $P\in \bdw$ and $\xi=\xi^P$ denote this cone by $\CC^P$. Then 
$$\CC^P\sbs Q^P(r_0,\rho_0) \q -\CC^P\sbs T^P(r_0,\rho_0)\sms \bar \Gw.$$
In a bounded Lipschitz domain, 
\begin{equation}\label{Phi0}
\rec{c_0}\gd(x)^{s_1}\leq \Phi_0(x)\leq c_0\gd(x)^{s_2},
\end{equation}
where $0<s_2\leq s_1$ depend on $\CC$.

Given $y\in T^P(r_0,\gr_0)\cap\bdw$ denote by $\mbfn_y^P$ the unit vector at $y$ pointing in the direction of the $\gx_1$ axis. We say that $\mbfn_y^P$ is the  \emph{approximate normal at $y$} associated with the local coordinates $\gx^P$.

\noindent\textit{Notation.}\hskip 2mm  Let $f_i$, $i=1,2$, be positive functions on some domain $X$. Then the notation
$f_1\sim f_2$ in $X$
means: there exists  $C>0$ \sth
$$\rec{C}f_1\leq f_2\leq C f_1 \qtxt{in }X.$$
The notation $f_1\lesssim f_2$ means: there exists  $C>0$ \sth
$f_1\leq Cf_2$ in $X$. The constant $C$ will be called a \textit{similarity constant}.

We state below the main results of this paper. Recall that $W:=\Phi_V/\Phi_0$.

 In the first theorem we present an estimate of $L_V$ harmonic functions.

\begin{theorem}\label{Knu-loc}  Assume (A1) and (A2). In addition assume that, for some $t^*>0$,

	\[J^* :=\int_{\Gw \cap [\gd<t^*]}\frac{\Phi_V^2}{\Phi_0}dx <\infty   \tag{B}  \]	

	Then there exists $C$ depending on $\bar a $, $\Gw$ \sth for every positive $\nu\in \GTM(\bdw)$,
	\begin{equation}\label{KV-global}
	\rec{C}\norm{\nu}\leq\int_{\Gw}\frac{\Phi_V}{\Phi_0}
	K_V[\nu] dx\leq C(1+J^*)\norm{\nu}.
	\end{equation}
\end{theorem}

The next theorem provides an estimate from below for $L_V$ potentials.

\begin{theorem}\label{th:I1}
	Assume (A1),  (A2). 
	Suppose that, in \eqref{Phi0}, $s_2<2$. Then there exists a  constant $c>0$ depending on $\bar a$, $r_0$, $\kappa$ and $s_2$ \sth, 
	
	\begin{equation}\label{G-I1}
	\int_{\Gw} \Phi_V d\tau\leq c \int_{\Gw}W G_V[\tau]dx  \forevery \tau \in \GTM_+(\Gw;\Phi_V).
	\end{equation}
		\end{theorem}

Following is an estimate from above of $L_V$ potentials. The estimate is obtained under two alternative sets of assumptions.

\begin{theorem}\label{th:I2} Assume (A1), (A2). 
	
\textbf{I.}  Suppose that there exist $\ga,\ga^*$ positive \sth, for every $P\in \bdw$:
	\[\rec{c}\gd^\ga\leq\Phi_V\leq c \gd^{\ga^* } \qtxt{in }\;Q^P(r_0/2,\gr_0/2). \tag{B2}\]
	In addition suppose that $s_1<2(1+\ga^*-\ga)$ ($s_1$ as in \eqref{Phi0}).
	
	Then there exists a  constant $c>0$ \sth, 
		\begin{equation}\label{G-I2}
	\int_{\Gw}W G_V[\tau]dx \leq c\int_{\Gw} \Phi_V d\tau \forevery \tau \in \GTM_+(\Gw;\Phi_V).
	\end{equation}
	
\textbf{II.}  Suppose that there exist constants $C_1>1$, $C_2>0$ \sth, for every $P\in \bdw$: 
	
	\noindent \emph{(B2')} \hskip 2mm If $x,z\in Q^P(r_0/2,\gr_0/2)$ and the segment $\overline{(x,z)}$ is parallel to the $\gx_1$ axis then,
	\begin{equation}\label{B2'}
	C_1 \gd(x)\leq \gd(z)<r_0/2 \Longrightarrow \Phi_V(x) \leq C_2\Phi_V(z).
	\end{equation}	
	Under these assumptions, if $s_1<2$ then \eqref{G-I2} holds.
	
	The constant in I. (respectively II.) depends on $\bar a$, $r_0$, $\kappa$, $s_1$ and the parameters in (B2) (respectively (B2')).
\end{theorem}

It is known that, under conditions (A1) and (A2), the following statement holds: 

Let  $u$ be either (i) a positive \LVsup function or (ii) a positive \LVsub function dominated by an \LVsup function. Then  $\tau:=- L_V u\in \GTM(\Gw;\Phi_V)$ and there exists a measure  $\nu\in \GTM(\bdw)$ \sth $u=K_V[\nu]+ G_V[\tau]$ (see Proposition 3.6 of \cite{MM2019N} which is valid in Lipschitz domains). The measure $\nu$ is called the $L_V$ boundary trace of $u$.

In view of the above statement, the following estimates are a direct consequence of the previous theorems.

\begin{theorem}\label{sub_sup} Assume (A1), (A2) and (B1), (B2').
In addition assume that $\Phi_0(x)\geq c_1 \gd(x)^{s_1}$ for some $s_1 \in (0,2)$.
	 Then there exists a constant $C>0$ depending  only on $\bar a$, $\Gw$ and the constants in (B2') \sth the following statements hold.
(a)	 Let $u$ be a positive \LVsup function, $\tau:=- L_V u$ and $\nu=\tr_Vu$. 
(see \eqref{Phi0}). Then,
\begin{equation}\label{sup_est}
\rec{C}(\int_\Gw \Phi_V d\tau +\norm{\nu})\leq \int_\Gw\,uWdx\leq C(\int_\Gw \Phi_V d\tau +\norm{\nu}).
\end{equation}

(b)	Let $u$ be a  positive \LVsub function dominated by an \LVsup function, $\tau:=- L_V u$ and $\nu=\tr_Vu$. (Note that $\tau$ is negative in this case.) Then,
\begin{equation}\label{sub_est}
\rec{C}\norm{\nu}\leq \int_\Gw\,uWdx + \int_\Gw \Phi_V d|\tau|\leq C\norm{\nu}.
\end{equation}	
\end{theorem}	

A similar result holds if $(B2')$ is replaced by (B2).

\section{Proof of Theorem \ref{Knu-loc}}

The proof is based on the following local result.

\begin{proposition}\label{p:Knu} Assume (A1), (A2).
	Let $P\in \bdw$,  let $\gx$ be a standard set of coordinates at $P$ and denote
	\begin{equation}
	J^P_\gb:=  \int_{\Gs^P_\gb(r_0,\gr_0)}\frac{\Phi_V^2}{\Phi_0}dS.
	\end{equation}
	
	Then there exists $t_1>0$ \sth, for every measure $\nu\in \GTM_+(\bdw)$ \sth $\supp \nu\sbs \bdw\cap T^P(r_0/4,\gr_0)$,
	\begin{equation}\label{e:KVy1}
	\rec{C}\norm{\nu}\leq\int_{\Gs^P_\gb(\frac{r_0}{2},\frac{\gr_0}{2})}\frac{\Phi_V}{\Phi_0}
	K_V[\nu] dS\leq C(1+J^P_\gb)\norm{\nu})  \forevery \gb\in (0,t_1)
	\end{equation}
	where $C$ and $t_1$ depend on $\bar a, \kappa, r_0$ but are independent of $P$.
\end{proposition}

\proof It is sufficient to show that for every $y\in \bdw\cap  T_P(r_0/4,\gr_0)$
	\begin{equation}\label{e:KVy2}
\rec{C}\leq\int_{\Gs^P_\gb(r_0/2,\gr_0/2)}\frac{\Phi_V}{\Phi_0}
K_V(\gx,y)dS\leq C (1+ J^P_\gb).
\end{equation}

It is known \cite {An87} that there exists $0< t_0<\gr_0$ \sth
\begin{equation}\label{e:Ky2}
K_V(y+t\mbfn^P_y, y) G_V(y+t\mbfn^P_y, x_0)\sim \, t^{2-N}  \q t\in (0,t_0).
\end{equation}
The similarity constant and $t_0$ depend on $\bar a, \Gw$ but not on $y$.\footnote{
This relation can also be derived from the estimates of the Green and Martin kernels in \cite{MM-Green}.}

Denote $\eta_y:=\gx^P-\gx^P(y)$ and
$$\CC_{y}^P :=\{\gx\in Q^P(r_0,\gr_0):4\gk|\eta'_y|<\eta_{y,1}\}.$$
If $t_1:=\min(t_0,\gk r_0)$ then, for every $y\in \bdw\cap  T_P(r_0/4,\gr_0)$ 
\begin{equation}\label{CC1}
\CC_{y,t_1}^P:=\CC_y^P \cap [\eta_{y,1}  < t_1] \sbs Q_{t_1}^P(r_0/2, \gr_0/2). 
\end{equation}

Using \eqref{e:Ky2}, \eqref{G-Phi} and the strong Hardy inequality we obtain:
\begin{equation}\label{e:Ky3}
\frac{\Phi_V}{\Phi_0}(\gx)\sim \frac{G_V}{G_0}(\gx,x_0)\sim \frac{K_0}{K_V}(\gx,y)  \forevery \gx\in \CC_{y,t_1}^P.
\end{equation}
Hence,
\begin{equation}\label{e:Ky4}
\frac{\Phi_V}{\Phi_0}(\gx)K_V(\gx,y)\sim K_0(\gx,y)  \forevery \gx\in \CC_{y,t_1}^P.
\end{equation}
In view of \eqref{Gsy} and the fact that the family of surfaces $\{\Gs^P_\gb(r_0,\gr_0)\}$ is uniformly Lipschitz, it follows that there exists a constant $C_1$ independent of $y$ and $\gb$ \sth
$$\rec{C_1}\leq\int_{\CC_{y,t_1}^P\cap \Gs^P_\gb(r_0,\gr_0)}\frac{\Phi_V}{\Phi_0}(\gx)
K_V(\gx,y)dS\leq C_1 \forevery \gb\in (0,t_1).$$

By BHP,
$$K_V(\cdot,y)\sim G_V(\cdot,x_0)\sim \Phi_V(\cdot) \qtxt{in } Q^P(r_0/2,\gr_0/2)\sms \CC_{y}^P  $$
Therefore, 
$$ \int_{\Gs^P_\gb(r_0/2,\gr_0/2)\sms \CC_y^P}\frac{\Phi_V}{\Phi_0}K_V(\gx,y)dS\leq J^P_\gb
\forevery \gb\in (0,t_1). $$
These relations imply \eqref{e:KVy2} which, by Fubini's theorem, implies \eqref{e:KVy1}.

\qed

\noindent\textit{Proof ot Theorem \ref{Knu-loc}.}\hskip 2mm  Let $P\in \bdw$ and let $\nu\in \GTM_+(\bdw)$ be a measure \sth $\supp\nu\sbs \bdw\cap T_P(r_0/4,\gr_0)$. Put $t_2:=min(t^*,t_1)$.

Then, integrating the terms in \eqref{e:KVy1} over $\gb\in (0,t_2)$ and using assumption (B) we obtain,
	\begin{equation}\label{KV-global.1}
\rec{C}t_2\norm{\nu}\leq\int_{Q^P_{t_2}(r_0/2,\gr_0/2)}\frac{\Phi_V}{\Phi_0}
K_V[\nu] dx\leq C(t_2+J^*)\norm{\nu}.
\end{equation}
	
For arbitrary $\nu\in \GTM_+(\bdw)$:
let $\{P_i\}_1^m$ be a set of points on $\bdw$ \sth $\{T^{P_i}(r_0/4,\gr_0/4)\}_1^m$ is a cover of $\bar \Gw\cap [0\leq \gd\leq t_2]$. Using a partition of unity for this cover, a measure $\nu$ can be written as a sum of measures $\nu_i$ supported in $T^{P_i}(r_0/4,\gr_0/4)\cap\bdw$, $i=1,\ldots,m$. Summing up inequalities \eqref{KV-global.1} for $P=P_i$, $i=1,\ldots,m$ we conclude that 
	\begin{equation}\label{KV-global.2}
\rec{C'}\norm{\nu}\leq\int_{A}\frac{\Phi_V}{\Phi_0} 
K_V[\nu] dx\leq C'(1+J^*)\norm{\nu}, 
\end{equation}
where $A:= \bigcup_1^m Q^{P_i}_{t_2}(r_0/2,\gr_0/2)$. Since $\Gw\cap [\gd<t_2/\bar c_\gk] \sbs A$, \eqref{Gsy} implies that $\Gw\sms A\sbs \Gw\cap[\gd\geq t_2/\bar c_\gk$.
Therefore
 $\Phi_V/\Phi_0$ and $K_V(\cdot,y)$ are bounded in $\Gw\sms A$. The bound depends on $\Gw$ and $\bar a$ but is independent of $V$. (This follows from the Harnack inequality and the fact that $\Phi_V(x_0)=K_V(x_0,y)=1$.)  Thus,
	\begin{equation}\label{KV-global.3} 
\int_{\Gw\sms A}\frac{\Phi_V}{\Phi_0} 
K_V[\nu] dx\leq C''\norm{\nu}.
\end{equation}
The last two inequalities imply \eqref{KV-global}. \qed

\vskip 4mm

\section{Proof of Theorem \ref{th:I1}}

The main part of the proof is contained in the following lemma which presents a local version of the global estimate.

\begin{lemma}\label{L-I1} Assume that (A1) and (A2) hold. Let $P\in \bdw$ and let $\gx=\gx^P$ be a standard set of coordinates in $T^P(r_0,\gr_0)$.
	Let $\tau$ be a positive measure in $\GTM_+(\Gw;\Phi_V)$ \sth $\supp\tau \sbs Q^P(r_0/4,\gr_0/4)$.
	Denote
		\begin{equation}\label{e:I1}
I_1^P(\gb):=\int_{\Gs^P_\gb(r_0/2,\gr_0/2)}\Phi_V(x)\int_\Gw G_V(x,y) \gc_{a\gb}(|x-y|)d\tau(y) dS_x,
\end{equation}
where $\gc_{r}(t)=\chr{(0,r)}(t)$ and $a\geq \max(16\kappa^2, 2\bar c_\gk)$.	Then there exists a constant $c$ depending on $a$, $\bar a$, $\bar c_\gk$ and $\Gw$ but independent of $P$ \sth,
	\begin{equation}\label{e:I1+}
	\rec{c}\int_{Q^P_{3a\gb/2}(r_0/4, \gr_0/4)}\Phi_V d\tau \leq \rec{\gb} I_1^P(\gb)\leq c\int_\Gw \Phi_V d\tau \forevery \gb\in (0, r_0/3a).
	\end{equation}
\end{lemma}

\proof 
We partition the domain $Q^P(r_0/2,\gr_0/2)$ into three parts $A_i$, $i=1,2,3$ and estimate each of the corresponding integrals separately. The sets are defined as follows:
$$\BAL
A_1=(Q^P_{\gb}\sms Q^P_{\gb/a^2})&(r_0/2,\gr_0/2),\q  A_2= Q^P_{\gb/a^2}(r_0/2,\gr_0/2),\\
&A_3=  (Q^P\sms Q^P_\gb)(r_0/2,\gr_0/2).
\EAL$$
Denote,
\begin{equation*}\BAL
I_{1,1}(\gb)&:= \int_{\Gs^P_\gb(r_0/2,\gr_0/2)}\Phi_V(x)\int_{A_1} G_V(x,y) \gc_{a\gb}(|x-y|)d\tau(y) dS_x,\\
I_{1,2}(\gb) &:= \int_{\Gs^P_\gb(r_0/2,\gr_0/2)}\Phi_V(x)\int_{A_2} G_V(x,y)
\gc_{a\gb}(|x-y|)d\tau(y) dS_x,\\
I_{1,3}(\gb)&:= \int_{\Gs^P_\gb(r_0/2,\gr_0/2)}\Phi_V(x)\int_{A_3} G_V(x,y) \gc_{a\gb}(|x-y|)d\tau(y)dS_x.
\EAL\end{equation*}
These integrals are estimated by essentially the same arguments as in the proof of 
\cite[Lemma 5.1]{MM2019N}. We present a sketch of the proof. In this inequality \eqref{Gsy} plays a crucial role.

Throughout the proof  $x\in \Gs^P_\gb$ and $|x-y|<a\gb$. 
The constants $c_i, c'_i, C_i$ etc. depend on $a$, $\bar c_\gk$, $\bar a$, $r_0$, $\gk$. By \eqref{Gsy},
\begin{equation}\label{gdx}
\gb/\bar c_\gk \leq \gd(x)\leq \gb, \q |x-y|<a\bar c_\gk \gd(x).
\end{equation}

If $y\in A_1$ then, by \eqref{Gsy}, $\gb/a^2\bar c_\gk\leq \gd(y)\leq \gb$. Therefore by the  Hardy (chain) inequality (see e.g. \cite[Lemma 3.2]{MM-Green})
\begin{equation}\label{Harnack1}
\rec{C_1}\Phi_V(x)\leq \Phi_V(y)\leq C_1\Phi_V(x). 
\end{equation}
In addition, by \cite[Theorem 1.3]{MM-Green},
\begin{equation}\label{Green1}
\rec{c_1}|x-y|^{2-N}\le G_V(x,y)\le c_1 |x-y|^{2-N}.
\end{equation}
Combining these inequalities and applying Fubini's theorem we obtain

\[\BAL  I_{1,1}(\gb)&\sim  \int_{\Gs^P_\gb(r_0/2,\gr_0/2)}\int_{A_1} |x-y|^{2-N}\gc_{a\gb}(|x-y|) \Phi_V(y) d\tau(y)dS_x \\
&=\int_{A_1} \int_{\Gs^P_\gb(r_0/2,\gr_0/2)}|x-y|^{2-N}\gc_{a\gb}(|x-y|)dS_x \Phi_V(y) d\tau(y). \EAL\]  
By an elementary computation, there exists a constant $c_1>0$ depending on $\gk$  (the Lipschitz constant of the surface $\Gs^P_\gb$) \sth 
$$\BAL  \rec{c_1} \gb  \leq \rec{c_1}(\frac{a}{\bar c_\gk} -1)\gb & \leq \int_{\Gs^P_\gb(r_0/2,\gr_0/2)}|x-y|^{2-N}\gc_{a\gb}(|x-y|)dS_x\\
&\leq c_1a\gb \forevery y\in A_1.
\EAL$$

Therefore
\begin{equation}\label{est-I_11}
I_{1,1}(\gb)\sim \gb \int_{A_1}\Phi_V d\tau
\end{equation}
with similarity constant depending on $a$, $\bar c_\gk$, $\bar a$,  $r_0$, $\gk$.

If  $y\in A_3$ then inequalities \eqref{Harnack1} and \eqref{Green1} remain valid. Therefore, in the same way as above, we obtain
\begin{equation}\label{est-I_13}
 I_{1,3}(\gb)\leq   c_3 \gb\int_{A_3}\Phi_V d\tau.
\end{equation}
We don't have a similar inequality from below for \emph{every} $y\in A_3$. If $y\in A_3$ is too far from $\Gs^P_\gb$  then $B_{a\gb}(y)\cap \Gs^P_\gb(r_0/2,\gr_0/2)$ may be empty. Of course we have an inequality for points that are close enough to $\Gs^P_\gb$. For instance if 
$A'_3 :=  (Q^P_{3a\gb/2}\sms Q^P_\gb)(r_0/2,\gr_0/2)$ then 
there exists a constant $c'>0$ \sth, for every $y\in A'_3$,
$$  B_{c'\gb}(y^*) \cap \Gs^P_\gb(r_0/2,\gr_0/2)   \sbs B_{a\gb}(y)$$
where $y^*$ is the projection of $y$ onto $\Gs^P_\gb(r_0/2,\gr_0/2)$  in the direction of the $\gx^P_1$ axis. Therefore there exists a constant $c'_3(a,\gk)>0$ \sth
$$\int_{\Gs^P_\gb(r_0/2,\gr_0/2)}|x-y|^{2-N}\gc_{a\gb}(|x-y|)dS_x>c'_3\gb \forevery y\in A'_3$$
Consequently 
\begin{equation}\label{est-I'13}
c'_3\gb\int_{A'_3}\Phi_V d\tau \leq I_{1,3}(\gb)   
\end{equation}

If $y\in A_2$ then, by \eqref{Gsy}, $\gd(y)<\gb/a^2$. Therefore, when $x\in \Gs^P_\gb$,
$$|x-y|>\gd(x)-\gd(y)>(\gb/\bar c_\gk)-\gd(y)>((a^2/\bar c_\gk)-1)\gd(y)\geq (2a-1) \gd(y).$$
Hence, by \cite[Theorem 1.4]{MM-Green},
\begin{equation}\label{est.G3}
\Phi_V(x)G_V(x,y) \sim \frac{\Phi_V(x)^2}{\Phi_V(x_y)^2}\Phi_V(y)|x-y|^{2-N},
\end{equation}
where $x_y\in \Gw$ may be chosen as follows:  the segment $\overline{(y,x_y)}$ is in the direction of the $\gx_1^P$ axis and $|x_y-y|=|x-y|$. Then 
$$|x_y-x|\leq  2|x-y|\leq 2a\gb, \; \gd(x)>\gb/\bar c_\gk \; \gd(y)> \gb/a^2\bar c_\gk.$$
Hence, by the strong Hardy inequality, there exists $c'(a,\gk)>0$ \sth
$$\rec{c'}\Phi_V(x)\leq \Phi_V(x_y)\leq c'\,\Phi_V(x).$$
Therefore, by \eqref{est.G3}, 
$$\Phi_V(x)G_V(x,y) \sim \Phi_V(y)|x-y|^{2-N},$$
for every $y\in A_2$ and $x\in \Gs^P_\gb(r_0/2,\gr_0/2)$ \sth $|x-y|<a\gb$.
Consequently,
\begin{equation}\label{est-I_12}\BAL
I_{1,2}(\gb)&\sim \int_{A_2} \int_{\Gs^P_\gb(r_0/2,\gr_0/2)}|x-y|^{2-N}\gc_{a\gb}(|x-y|)dS_x \Phi_V(y) d\tau(y)\\
&\sim \gb\int_{A_2}\Phi_V d\tau.
\EAL\end{equation}

Since $Q^P_{3a\gb/2}\gb)(r_0/2,\gr_0/2) = A_1\cup A_2\cup A'_3$, the left hand inequality in \eqref{e:I1+}
follows from \eqref{est-I_11}, \eqref{est-I'13} and \eqref{est-I_12}. The right hand inequality in \eqref{e:I1+} follows from  \eqref{est-I_11}, \eqref{est-I_13} and \eqref{est-I_12}.

\qed

\noindent\textit{Proof of Theorem \ref{th:I1}.} \hskip 2mm

By Lemma \ref{L-I1}, 
$$	\rec{c}\int_{Q^P_{3a\gb/2}(\frac{r_0}{4}, \frac{\gr_0}{4})}\Phi_V d\tau \leq I_1(\gb)$$
for every $\gb<r_0/3a$. Therefore for every $r_0/6a<\gb<r_0/3a$,
$$	(\gb/c)\int_{Q^P_{\frac{r_0}{4}}(\frac{r_0}{4}, \frac{\gr_0}{4})}\Phi_V d\tau \leq I^P_1(\gb) \forevery \gb\in (r_0/6a, r_0/3a).$$
Hence, by \eqref{Phi0},
$$ \rec{c}\gb^{1-s_2} \int_{Q^P_{\frac{r_0}{4}}(\frac{r_0}{4}, \frac{\gr_0}{4})}\Phi_V d\tau \leq \int_{\Gs^P_\gb(r_0/2,\gr_0/2)}\frac{\Phi_V(x)}{\Phi_0(x)}\int_\Gw G_V(x,y) d\tau(y) dS_x.$$
Since, by assumption, $s_2<2$, integration over $\gb$ in the interval $(r_0/6a, r_0/3a)$ yields:
$$c_1 \int_{Q^P_{\frac{r_0}{4}}(\frac{r_0}{4}, \frac{\gr_0}{4})}\Phi_V d\tau \leq \int_{\Gs^P_\gb(r_0/2,\gr_0/2)}\frac{\Phi_V(x)}{\Phi_0(x)}\int_\Gw G_V(x,y) d\tau(y) dS_x,$$
where $c_1$ depends on $s_2$ and the choice of $a$. Choosing for instance $a=2\max(16\gk^2, C_1)$ we obtain 
	\begin{equation}\label{G-I1'}
\int_{Q^P_{\frac{r_0}{4}}(\frac{r_0}{4}, \frac{\gr_0}{4})} \Phi_V d\tau\leq c' \int_{\Gw}W G_V[\tau]dx  \forevery \tau \in \GTM_+(\Gw;\Phi_V).
\end{equation}
 where the constant depends on $r_0$, $\gk$, $s_2$, $\bar a$ but is indepent of $P$.
 
Let $\{P_i\}_1^m$ be a set of points on $\bdw$ \sth 
 $$\Gw_{r_0/4\bar c_\gk}:=\{x\in \Gw:  \gd(x)< r_0/4\bar c_\gk   \}\sbs \cup_1^m \{Q^{P_i}_{\frac{r_0}{4}}(r_0/4,\gr_0)\}$$ (see \eqref{Phi0}). Using a partition of unity corresponding to this cover, a positive measure $\tau\in \GTM(\Gw;\Phi_V)$ can be written as a sum 
$\tau= \sum_1^m\tau_i+\tau'$ \sth 
$$\supp \tau_i \sbs Q^P_{\frac{r_0}{4}}(\frac{r_0}{4}, \frac{\gr_0}{4}),\q \supp\tau'\sbs \{y\in \Gw: \gd(y)\geq   r_0/4\bar c_\gk  \}. $$ 
	Summing up inequalities \eqref{G-I1'} for $P=P_i$,
	$i=1\ldots,m$ we obtain
\begin{equation}
\int_{\Gw_{r_0/4\bar c_\gk}} \Phi_V d\tau\leq c' \int_{\Gw}W G_V[\tau]dx  \forevery \tau \in \GTM_+(\Gw;\Phi_V).
\end{equation}

On the other hand,
\begin{equation}\label{temp3}\BAL
&\int_{\Gw}W(x) \int_\Gw G_V(x,y)d\tau'(y) dx =
\int_{\Gw} \int_\Gw \frac{\Phi_V}{\Phi_0}(x) G_V(x,y) dx d\tau'(y)\geq\\
&\int_{[\gd(y)\geq r_0/4\bar c_\gk]} \int_{[|x-y|<r_0/8\bar c_\gk]} \frac{\Phi_V}{\Phi_0}(x) G_V(x,y) dx d\tau'(y).
\EAL\end{equation}
 If $|x-y|<r_0/8\bar c_\gk$ and $\gd(y)>r_0/4\bar c_\gk$ then, by the strong Harnack inequality, $$\rec{c_0}\Phi_V(x)\leq \Phi_V(y)\leq c_0\Phi_V(x)$$ 
 and by \cite[Theorem 1.3]{MM-Green}, $\G_V(x,y)\leq c |x-y|^{2-N}$. 
 Therefore by \eqref{temp3},
\begin{equation*}\BAL
&\int_{\Gw}\frac{\Phi_V}{\Phi_0}(x) \int_\Gw G_V(x,y)d\tau'(y) dx \geq\\
\frac{c_1}{r_0^{s_2}}&\int_{[\gd(y)\geq r_0/4\bar c_\gk]} \int_{|x-y|<r_0/8} |x-y|^{2-N} dx \,\Phi_V(y)d\tau'(y)\geq
c_2  \int_{\Gw}\Phi_V d\tau',
\EAL \end{equation*}
This inequality and \eqref{temp3} imply \eqref{G-I1}.
\qed

\vskip 3mm

\section{Proof of Theorem \ref{th:I2}}

The proof is based on Lemma \ref{L-I1} and the next two lemmas

\begin{lemma}\label{L-I2}
	Assume (A1), (A2). Let $P\in \bdw$ and assume that there exist $\ga,\ga^*$ positive \sth (B2) holds.
Without loss of generality we assume that $\ga\neq 1/2$. (If $\ga=1/2$ we replace it by a slightly larger exponent.)

Let
	\begin{equation}\label{e:I2}
I^P_2(\gb):=	\int_{\Gs^P_\gb(r_0/2,\gr_0/2)}\Phi_V(x)\int_\Gw G_V(x,y) (1-\gc_{a\gb}(|x-y|))d\tau(y) dS_x
	\end{equation}
	with $a$ as in  Lemma \ref{L-I1}.
	Then there exists $C>0$  \sth,  if $\tau\in \GTM_+(\Gw;\Phi_V)$, $\supp \tau \sbs Q^P(r_0/4,\gr_0/4)$ and $\gb\in (0,  r_0/4)$ then,
	\begin{equation}\label{Gtau-loc}\BAL
I^P_2(\gb) 
\leq C\gb^\gl\int_{\Gw}\Phi_V d\tau,
	\EAL	\end{equation}
	where $\gl := 2(\ga^*-\ga)+1$.
\end{lemma}

\proof
In the domain of integration of $I^P_2$, $\gd(x)\leq \gb$ and $|x-y|\geq a\gb$.
Therefore estimate \eqref{est.G3} holds with $x_y$ as follows:  the segment $\overline{(x,x_y)}$ is in the direction of the $\gx_1^P$ axis and $|x_y-x|=|x-y|$.

By assumption (B2),
$$\frac{\Phi_V(x)}{\Phi_V(x_y)} \leq c(a)\frac{\gb^{\ga^*}}{|x-y|^{\ga}}.$$
Hence, by   \eqref{est.G3}, 

\begin{equation}\label{Gtau-Lip}\BAL
&I^P_2(\gb)=\\
&\int_\Gw \int_{\Gs^P_\gb(\frac{r_0}{2},\frac{\gr_0}{2})}\frac{\Phi_V(x)^2}{\Phi_V(x_y)^2}|x-y|^{2-N} (1-\gc_{a\gb}(|\xi-\eta|))dS_x \Phi_V(y)d\tau(y)\lesssim\\
&\int_\Gw \int_{\Gs^P_\gb(\frac{r_0}{2},\frac{\gr_0}{2})}\gb^{2\ga^*}|x-y|^{2-N-2\ga} (1-\gc_{a\gb}(|\xi-\eta|))dS_x \Phi_V(y)d\tau(y)\lesssim\\
&\int_\Gw \big(\int_{a\gb}^{\diam(\Gw)}r^{-2\ga}dr\big) \gb^{2\ga^*}\Phi_V(y)d\tau(y)\lesssim
\gb^{\gl}\int_{\Gw}\Phi_V d\tau.
\EAL	\end{equation}

\qed

A similar  estimate of $I^P_2(\gb)$  can be obtained if condition (B2) is replaced by (B2').

\begin{lemma}\label{L-I2'}
	Assume (A1), (A2). Let $P\in \bdw$ and assume that (B2') holds. 	

	Let $I^P_2(\gb)$ be defined as in \eqref{e:I2} with  $a\geq \max(16\kappa^2, 2\bar c_\gk, C_1)$. 
	Then there exists $C>0$ \sth  if $\tau\in \GTM_+(\Gw;\Phi_V)$, $\supp \tau \sbs Q^P(r_0/4,\gr_0/4)$ and $\gb\in (0,  r_0/4)$,
	\begin{equation}\label{Gtau-loc}\BAL
	I^P_2(\gb) 
	\leq C\gb\int_{\Gw}\Phi_V d\tau.
	\EAL	\end{equation}
\end{lemma}

\proof Estimate \eqref{est.G3} holds with $x_y$ as before.  In the domain of integration,
$$
|x_y-x|= |x-y| \geq a\gb \geq C_1\gd(x) \forevery x\in \Gs^P_\gb(r_0/2,\gr_0/2)$$
Therefore by \eqref{est.G3} and \eqref{B2'}
$$\BAL
I^P_2(\gb) &\sim \int_\Gw \int_{\Gs^P_\gb(\frac{r_0}{2},\frac{\gr_0}{2})}|x-y|^{2-N} (1-\gc_{a\gb}(|x-y))dS_x \Phi_V(y)d\tau(y)\\
& \lesssim \gb \int_\Gw \Phi_V d\tau.
\EAL$$

\qed

\noindent\textit{Proof of Theorem \ref{th:I2}.} First we assume that $\supp \tau \sbs Q^P(r_0/4,\gr_0/4)$ for some $P\in \bdw$.  
Recall that by Lemma  \ref{L-I1}, inequality \eqref{e:I1+}, we have
$$\BAL
I_1^P(\gb):=&\int_{\Gs^P_\gb(r_0/2,\gr_0/2)}\Phi_V(x)\int_\Gw G_V(x,y) \gc_{a\gb}(|x-y|)d\tau(y) dS_x\\
\leq&
c\gb\int_\Gw\Phi_V d\tau
\EAL$$
for every $a\geq a_\gk:=\max(16\kappa^2, 2\bar c_\gk)$ and $\gb\in (0, r_0/3a)$. Hence, by \eqref{Gsy} and \eqref{Phi0},
\begin{equation}\label{temp5}
\BAL
&\int_{\Gs^P_\gb(r_0/2,\gr_0/2)}\frac{\Phi_V}{\Phi_0}(x)\int_\Gw G_V(x,y) \gc_{a\gb}(|x-y|)d\tau(y) dS_x\\
&\leq c\gb^{1-s_1}\int_\Gw\Phi_V d\tau \forevery \gb\in (0,r_\gk)
\EAL\end{equation}
where $r_\gk:= r_0/3a_\gk$.

 If (B2) holds then, by Lemma \ref{L-I2},
$$	\int_{\Gs^P_\gb(r_0/2,\gr_0/2)}\frac{\Phi_V}{\Phi_0}\int_\Gw G_V(x,y) (1-\gc_{a\gb}(|x-y|))d\tau(y) dS_x \lesssim \gb^{\gl - s_1}\int_\Gw \Phi_V\,d\tau $$
with $a$ as in Lemma \ref{L-I1}. In particular, choosing $a=a_\gk$, 
this inequality  and \eqref{temp5} yield, (as $\gl\leq 1$)
\begin{equation}
\int_{\Gs^P_\gb(r_0/2,\gr_0/2)}W \int_\Gw G_V(x,y)d\tau(y) dS_x \lesssim \gb^{\gl - s_1}\int_\Gw \Phi_V\,d\tau 
\end{equation}
for every $\gb\in (0,r_\gk)$. By assumption $s_1<2(1+\ga^*-\ga)= 1+\gl$. Thus $-1<\gl-s_1$ and 
integrating over $\gb$ in $(0,r_\gk)$ we obtain,
\begin{equation}\label{base1}
	\int_{Q^P_{r_\gk}(r_0/2,\gr_0/2)}W(x)\int_\Gw G_V(x,y) d\tau(y) dx \leq c\int_\Gw \Phi_V\,d\tau 
\end{equation}
where  the constant dedepends on $\bar a$, $r_0$, $\gk$, $\gl$ and $s_1$. 

Similarly, if (B2') holds then, by Lemma \ref{L-I2'},
$$	\int_{\Gs^P_\gb(r_0/2,\gr_0/2)}\frac{\Phi_V}{\Phi_0}\int_\Gw G_V(x,y) (1-\gc_{a\gb}(|x-y|))d\tau(y) dS_x \lesssim \gb^{1 - s_1}\int_\Gw \Phi_V\,d\tau.$$
Combining this inequality with \eqref{temp5} we obtain,
\begin{equation}
\int_{\Gs^P_\gb(r_0/2,\gr_0/2)}W \int_\Gw G_V(x,y)d\tau(y) dS_x \lesssim \gb^{1 - s_1}\int_\Gw \Phi_V\,d\tau 
\end{equation}
In this case we  assume that $-2<s_1$. Therefore integrating over $\gb$ in $(0,r_\gk)$ we again obtain \eqref{base1}.

Next we show that, 
\begin{equation}\label{comp1}
\int_{D^P_{r_\gk}(r_0/2,\gr_0/2)}W(x)\int_\Gw G_V(x,y) d\tau(y) dx \leq c\int_\Gw \Phi_V\,d\tau 
\end{equation}
where $D^P_{t}(r,\gr):= Q^P(r,\gr) \sms Q^P_t(r,\gr)$. 

By \eqref{Gsy}, if $x\in D^P_{r_\gk}(r_0/2,\gr_0/2)$ then 
 $\gd(x)\geq r'_\gk:=r_\gk/\bar c_\gk$. Therefore, by Harnack's inequality, $W$ is bounded and bounded away from zero in this set, by  constants dependent on $r'_\gk$ but independent of $V$. (Recall that $\Phi_V(x_0)=1$.) 

Denote, 
$$E := Q^P_{r'_\gk/2}(r_0/2,\gr_0/2), \;E':= Q^P(r_0/2,\gr_0/2) \sms E = D^P_{r'_\gk/2}(r_0/2,\gr_0/2).$$
There exists a constant $C>0$ dependent on $r_\gk$, $\bar c_\gk$ and $\bar a$ \sth 
$$\rec{C}\Phi_V(y)\leq G_V(x,y)\leq C\Phi_V(y) \forevery (x,y)\in D^P_{r_\gk}(r_0/2,\gr_0/2)\ti E.$$
Therefore
 \begin{equation}\label{comp3}
\int_{D^P_{r_\gk}(r_0/2,\gr_0/2)}W(x)\int_E G_V(x,y) d\tau(y) dx \leq c_1\int_E \Phi_V\,d\tau. 
\end{equation}
On the other hand, since $E'\sbs \{x\in \Gw: \gd(x)>r'_\gk/2\bar c_\gk\}$, there exists a constant $C'$ depending on $r_\gk$, $\bar c_\gk$ and $\bar a$ \sth,
$$G_V(x,y)\leq C'|x-y|^{2-N} \forevery (x,y)\in (E')^2.$$
 Therefore
 \begin{equation}\label{comp4}
\int_{D^P_{r_\gk}(r_0/2,\gr_0/2)}W(x)\int_{E'} G_V(x,y) d\tau(y) dx \leq c'\int_{E'} \Phi_V\,d\tau. 
\end{equation}
Inequalities \eqref{comp3} and \eqref{comp4} yield \eqref{comp1}. This inequality and \eqref{base1} imply
\begin{equation}\label{temp6}
\int_{Q^P(r_0/2,\gr_0/2)} W(x) G_V[\tau]dx \leq c\int_\Gw \Phi_V d\tau
\end{equation}
 for measures $\tau\in \GTM_+(\Gw;\Phi_V)$ \sth $\supp \tau \sbs Q^P(r_0/4,\gr_0/4)$.

Finally, let $\tau$ be an arbitrary measure in $\GTM_+(\Gw;\Phi_V)$. Let $\{P_i\}_1^m$ be a subset of $\bdw$ \sth  $\cup_1^m  T^{P_i}(r_0/4,\gr_0/4)$ covers a set $A_t= \Gw\cap [\gd<t]$
for some $t>0$ depending on $r_0, \gk$.
 Using a corresponding partition of unity we write $\tau\chr{A}=\sum_1^m\tau_i$ where $\supp\tau_i\in Q^{P_i}(r_0/4, \gr_0/4)$. Summing up the terms in \eqref{temp6} over $P=P_i$, $i=1,\ldots, m$ we obtain 
\begin{equation}\label{temp7}
\int_{A} W(x) G_V[\tau]dx \leq c\int_\Gw \Phi_V d\tau \forevery \tau\in \GTM_+(\Gw;\Phi_V).
\end{equation}
As $W$ is bounded in $\Gw\sms A$,
$$\int_{\Gw \sms A} W(x) G_V[\tau]dx \leq c\int_\Gw \Phi_V d\tau \forevery \tau\in \GTM_+(\Gw;\Phi_V).$$
 This inequality and \eqref{temp7} imply \eqref{G-I2}.
\qed

\end{document}